# RBF Interpolation with CSRBF of Large Data Sets

Vaclav Skala
*University of West Bohemia, Plzen, Czech Republic.*
*www.VaclavSkala.eu*

**Abstract**
This contribution presents a new analysis of properties of the interpolation using Radial Bases Functions (RBF) related to large data sets interpolation. The RBF application is convenient method for scattered $d$-dimensional interpolation. The RBF methods lead to a solution of linear system of equations and computational complexity of solution is nearly independent of a dimensionality. However, the RBF methods are usually applied for small data sets with a small span geometric coordinates. This contribution explores properties of the RBF interpolation for large data sets and large span of geometric coordinates of the given data sets with regard to expectable numerical stability of computation.

*Keywords:* Radial basis functions, compactly supported RBF, interpolation, meshless interpolation.

## 1  Introduction

Interpolation techniques are used in solution of many engineering problems. However, the interpolation of unorganized scattered data is still a severe problem. The standard approaches are based on tessellation of the domain in $x, y$ or $x, y, z$ spaces using e.g. Delaunay triangulation etc. This approach is applicable for static data and $t$-varying data, if data in the time domain are "framed", i.e. given for specific time samples. However, it leads to increase of dimensionality, i.e. from triangulation in $E^2$ to triangulation in $E^3$ or from triangulation in $E^3$ to triangulation in $E^4$ etc. It results to significant increase of triangulation complexity and an algorithm implementation.

Radial Basis Functions (RBF) offer several significant advantages, e.g. RBF formulation leads to a solution of linear system of equations, i.e. $\boldsymbol{Ax} = \boldsymbol{b}$, and it is applicable generally to $d$-dimensional problems. It does not require tessellation of the definition domain and it is especially convenient for scattered data interpolation. The RBF interpolation is smooth by a definition and it can be applied for interpolation of scalar and vector fields. If the Compactly Supported RBFs (CSRBF) are used, sparse matrix data structures can be used which decreases memory requirements significantly.

However, there are some weak points of a RBF application in real problems solution, e.g. problems with robustness of the RBF due to low conditionality of the matrix $\boldsymbol{A}$ and with numerical stability, if the interpolation is to be applied over a large span of $x$ values. Memory requirements are of $O(N^2)$ complexity, where $N$ is a number of points in which values are given. Computational complexity of a solution of the LSE, which is $O(N^3)$, resp. $O(kN^2)$, where $k$ is a number of iteration if iterative method is used. Also severe problems with unexpected behavior at geometrical borders, e.g. in PDE solutions.

The meshless techniques are easily scalable to higher dimensions and can handle scattered spatial and spatial-temporal data. The polygonal based techniques require tessellations, e.g. Delaunay triangulation with $O\left(N^{\lfloor d/2+1 \rfloor}\right)$ computational complexity for $N$ points in $d$-dimensional space or another tessellation method. However, the complexity of tessellation algorithms implementation grow significantly with dimensionality and severe problems with robustness might be expected as well.

## 2  RBF Interpolation

Meshless (meshfree) methods are based on the idea of Radial Basis Function (RBF) interpolation [2], [22], [22], [16], which is not separable, and lead to a solution of a linear system equations (LSE) with a full or sparse matrix [4], [5]. The RBF interpolation is based on computing of the distance of two points in the $d$-dimensional space and is defined by a function:

$$f(x) = \sum_{j=1}^{N} \lambda_j \, \varphi(\|x - x_j\|) = \sum_{j=1}^{N} \lambda_j \, \varphi(r_j)$$

Formally: $r_j = \|x - x_j\|_2 \stackrel{\text{def}}{=} \sqrt{(x - x_j)^2 + (y - y_j)^2} = \sqrt{(x - x_j)^2 + (y - y_j)^2 + (1-1)^2}$

$x = [x, y : 1]^T$ is expressed in homogeneous coordinates for $d = 2$ case, $\lambda_j$ are weights to be computed. In the case of the projective space representation $r_{ij} = \|x_i - x_j\|_P \stackrel{\text{def}}{=} \|w_j x_i - w_i x_j\|_P = \sqrt{(w_j x_i - w_i x_j)^2 + (w_j y_i - w_i y_j)^2 + (w_i w_j - w_i w_j)^2} = \sqrt{(w_j x_i - w_i x_j)^2 + (w_j y_i - w_i y_j)^2}$

It means that for the given data set $\{\langle x_i, h_i \rangle\}_1^N$, where $h_i$ are associated values to be interpolated and $x_i$ are domain coordinates, we obtain a linear system of equations:

$$h_i = f(x_i) = \sum_{j=1}^{N} \lambda_j \, \varphi(\|x_i - x_j\|) + P_k(x_i) \qquad x = [x, y:1]^T \qquad i = 1, \dots, N$$

Due to some stability issues, usually a polynomial $P_k(x)$ of a degree $k$ is added. Mostly, a polynomial of the 1st degree is used, i.e. $P_1(x) = \boldsymbol{a}^T x$.

$$f(x_i) = \sum_{j=1}^{N} \lambda_j \, \varphi(\|x_i - x_j\|) + \boldsymbol{a}^T x_i = \sum_{j=1}^{N} \lambda_j \, \varphi_{i,j} + \boldsymbol{a}^T x_i \qquad h_i = f(x_i) \qquad i = 1, \dots, N$$

and additional conditions are to be applied:

$$\sum_{j=1}^{N} \lambda_i x_i = \boldsymbol{0} \qquad \text{i.e.} \qquad \sum_{j=1}^{N} \lambda_i x_i = 0 \qquad \sum_{j=1}^{N} \lambda_i y_i = 0 \qquad \sum_{j=1}^{N} \lambda_i = 0$$

It leads to s linear system of equations

$$\begin{bmatrix} \varphi_{1,1} & \cdots & \varphi_{1,N} & x_1 & y_1 & 1 \\ \vdots & \ddots & \vdots & \vdots & \vdots & \vdots \\ \varphi_{N,1} & \cdots & \varphi_{N,N} & x_N & y_N & 1 \\ x_1 & \cdots & x_N & 0 & 0 & 0 \\ y_1 & \cdots & y_N & 0 & 0 & 0 \\ 1 & \cdots & 1 & 0 & 0 & 0 \end{bmatrix} \begin{bmatrix} \lambda_1 \\ \vdots \\ \lambda_N \\ a_x \\ a_y \\ a_0 \end{bmatrix} = \begin{bmatrix} h_1 \\ \vdots \\ h_N \\ 0 \\ 0 \\ 0 \end{bmatrix} \qquad \begin{bmatrix} \boldsymbol{B} & \boldsymbol{P} \\ \boldsymbol{P}^T & \boldsymbol{0} \end{bmatrix} \begin{bmatrix} \boldsymbol{\lambda} \\ \boldsymbol{a} \end{bmatrix} = \begin{bmatrix} \boldsymbol{f} \\ \boldsymbol{0} \end{bmatrix}$$

$$Ax = b$$
$$\boldsymbol{a}^T x_i = a_x \, x_i + a_y \, y_i + a_0$$

It can be seen that for $d$-dimensional case a system of $(N + d + 1)$ LSE has to be solved, where $N$ is a number of points in the dataset and $d$ is the dimensionality of data. For $N$ points given and $d = 2$, i.e. $\boldsymbol{a} = [a_x, a_y: a_0]^T$, a system of $(N + 3)$ linear equations has to be solved. If "global" functions, e.g. TPS ($\varphi(r) = r^2 \lg r$, are used, then the matrix $\boldsymbol{B}$ is "full", if "local" functions (Compactly Supported RBF – CSRBF) are used, the matrix $\boldsymbol{B}$ can be sparse.

The RBF interpolation was originally introduced by multiquadric method in 1971 [5], which was called Radial Basis Function (RBF) method. Since then many different RFB interpolation schemes have been developed with some specific properties, e.g. Thin-Plate Spline (TPS uses $\varphi(r) = r^2 \lg r$) [4], a function $\varphi(r) = e^{-(\epsilon r)^2}$ was proposed in [22] and CSRBFs were introduced in [21] as:

$$\varphi(r) = \begin{cases} (1-r)^q P(r), & 0 \leq r \leq 1 \\ 0, & r > 1 \end{cases},$$

where: $P(r)$ is a polynomial function and $q$ is a parameter. Theoretical problems with numerical stability were solved in [4]. The CSRBFs are defined for the "normalized" interval $r \in \langle 0, 1 \rangle$, but for the practical use a scaling is used, i.e. the value $r$ is multiplied by a scaling factor $\alpha > 0$ (shape parameter).

Meshless techniques are used in engineering problem solutions, e.g. partial differential equations (PDE) [6] surface modeling [8], surface reconstruction of scanned objects [3], [18] reconstruction of corrupted images [22], etc. More generally, meshless object representation is based on specific interpolation techniques [1][2][6][19][22]. In addition, subdivision or hierarchical methods are used to decrease sizes of computations and increase robustness [14][20]. If "global" RBF functions are considered, the RBF matrix is full and in the case of $10^6$ of points, the RBF matrix is of the size approx.$10^6 \times 10^6$ ! On the other hand, if CSRBF used, the relevant matrix is sparse and computational and memory requirements can be decreased significantly using special data structures.

## 3  Decomposition of RBF Interpolation

The RBF interpolation can be described in the matrix form as

$$\begin{bmatrix} A & P \\ P^T & 0 \end{bmatrix} \begin{bmatrix} \lambda \\ a \end{bmatrix} = \begin{bmatrix} f \\ 0 \end{bmatrix} \qquad \boldsymbol{a}^T \boldsymbol{x}_i = a_x x_i + a_y y_i + a_0$$

where $\boldsymbol{x} = [x, y: 1]^T$, the matrix $\boldsymbol{A}$ is symmetrical and semidefinite positive (or strictly positive) definite. Let us consider the Schur's complement (validity of all operation is expected)

$$\boldsymbol{M} = \begin{bmatrix} A & B \\ C & D \end{bmatrix} = \begin{bmatrix} I & 0 \\ CA^{-1} & I \end{bmatrix} \begin{bmatrix} A & 0 \\ 0 & M/A \end{bmatrix} \begin{bmatrix} I & A^{-1}B \\ 0 & I \end{bmatrix} \qquad M/A \stackrel{\text{def}}{=} D - CA^{-1}B$$

where $\boldsymbol{M}/\boldsymbol{A}$ is the Schur's complement. Then the inversion matrix $\boldsymbol{M}^{-1}$

$$\boldsymbol{M}^{-1} = \begin{bmatrix} I & -A^{-1}B \\ 0 & I \end{bmatrix} \begin{bmatrix} A^{-1} & 0 \\ 0 & (M/A)^{-1} \end{bmatrix} \begin{bmatrix} I & 0 \\ -CA^{-1} & I \end{bmatrix}$$

Now, the Schur's complement can be applied to the RBF interpolation. As the matrix $\boldsymbol{M}$ is nonsingular, inversion of the matrix $\boldsymbol{M}$ can be used. Using the Schur's complement (as the matrix $\boldsymbol{D} = \boldsymbol{0}$)

$$\boldsymbol{M}^{-1} = \begin{bmatrix} I & -A^{-1}P \\ 0 & I \end{bmatrix} \begin{bmatrix} A^{-1} & 0 \\ 0 & (M/A)^{-1} \end{bmatrix} \begin{bmatrix} I & 0 \\ -P^TA^{-1} & I \end{bmatrix} \qquad M/A \stackrel{\text{def}}{=} P^T A^{-1} P$$

Then $\det(\boldsymbol{M}) \neq 0$, $\det(\boldsymbol{A}) \neq 0$ and $\det(\boldsymbol{M}/\boldsymbol{A}) \neq 0$ as the matrices are nonsingular.

However, if RBF interpolation is used for larger data sets, there is a severe problem with robustness and numerical stability, i.e. numerical computability issues. Using the Schur's complement we obtain

$$\det(\boldsymbol{M}^{-1}) = \frac{1}{\det(\boldsymbol{M})} = \frac{1}{\det(\boldsymbol{A})} \frac{1}{\det(\boldsymbol{M}/\boldsymbol{A})} = \frac{1}{\det(\boldsymbol{A})} \frac{1}{\det(\boldsymbol{P}^T \boldsymbol{A}^{-1} \boldsymbol{P})}$$

Properties of the matrix $\boldsymbol{A}$ are determined by the RFB function used. The value of $\det(\boldsymbol{A})$ depends also on the mutual distribution of points. However, the influence of $\det(\boldsymbol{P}^T\boldsymbol{A}^{-1}\boldsymbol{P})$ is also significant as the value depends on the points mutual distribution due to the matrix $\boldsymbol{A}$ but also to points distribution in space, due to the matrix $\boldsymbol{P}$. It means that translation of points in space does have significant influence. Let us imagine for a simplicity that the matrix $\boldsymbol{A} = \boldsymbol{I}$ (it can happen if CSRBF is used and only one point is within the radius $r = 1$). Then the distance of a point from the origin has actually quadratic influence, if a polynomial $P_1(x)$ is used, as the point's position is in the matrices $\boldsymbol{P}^T$ and $\boldsymbol{P}$. There is a direct significant consequence for the RBF interpolation.

$$f(x) = \sum_{j=1}^{N} \lambda_j \, \varphi(\|x - x_j\|) + P_k(x)$$

when the $P_k(x)$, $k = 1, 2$ is a polynomial

$$P_1(x) = a_0 + a_1 x + a_2 y \quad , \text{resp.} \quad P_2(x) = a_0 + a_1 x + a_2 y + a_3 x^2 + a_4 xy + a_5 y^2$$

In the case of $\boldsymbol{A} = \boldsymbol{I}$, we get a matrix $\boldsymbol{P}^T\boldsymbol{P}$ of the size $(3 \times 3)$ and $\det(\boldsymbol{P}^T\boldsymbol{P})$ in the case of a linear polynomial $P_1(x)$ as:

$$\det(\boldsymbol{P}^T\boldsymbol{P}) = \begin{vmatrix} \sum_{i=1}^{N} x_i^2 & \sum_{i=1}^{N} x_i y_i & \sum_{i=1}^{N} x_i \\ \sum_{i=1}^{N} x_i y_i & \sum_{i=1}^{N} y_i^2 & \sum_{i=1}^{N} y_i \\ \sum_{i=1}^{N} x_i & \sum_{i=1}^{N} y_i & \sum_{i=1}^{N} 1 \end{vmatrix} = n\left(\sum x_i^2 \sum y_i^2\right) - \sum y_i(\dots) + \sum y_i(\dots) \dots$$

It means that points distribution in space and their distances from the origin play a significant role as the $\det(\boldsymbol{P}^T\boldsymbol{P})$ contains elements $\sum_{i=1}^{N} x_i^2$ and $\sum_{i=1}^{N} y_i^2$ in multiplicative etc. in the linear polynomial case.

If a quadratic polynomial $P_2(x)$ is used, the matrix $\boldsymbol{P}^T\boldsymbol{P}$ is of the size $(6 \times 6)$ and $\det(\boldsymbol{P}^T\boldsymbol{P})$ contains elements $\sum_{i=1}^{N} x_i^4$, $\sum_{i=1}^{N} y_i^2$,..., $\sum_{i=1}^{N} 1$ in multiplicative, which brings even worst situation as the matrix $\boldsymbol{P}^T\boldsymbol{P}$ contains small and very high values.

As a direct consequence, eigenvalues will have large span and therefore the linear system of equations will become extremely ill-conditioned.

## 4 Conclusion

The RBF interpolation using compactly supported RBF (CSRBF) have several significant advantages over methods based on smooth interpolation made on triangulated space area. In this contribution, some properties of the CSRBF interpolation methods have been presented from "engineering" point of view and selected features related to robustness and stability of computation have been presented. The presented founding are fundamental especially in the case of engineering applications. In future, CSRBF analysis of sparse data structures used, space subdivision and speed up of computation will be explored more deeply especially for RBF based approximations.

## Acknowledgment

The author would like to thank students and colleagues at the University of West Bohemia, Plzen, their valuable comments and hints provided, to anonymous reviewers for their critical view that helped to improve the manuscript. This research was supported by the projects GACR project No.GA 17-05534S.

# References

[1] Adams,B., Ovsjanikov,M., Wand,M., Seidel,H.-P., Guibas,L.J.: Meshless modeling of Deformable Shapes and their Motion. ACM SIGGRAPH Symp. on Computer Animation, 2008

[2] Buhmann,M.D.: Radial Basis Functions:Theory&Implementations.Cambridge Univ.Press,2008.

[3] Carr,J.C., Beatsom,R.K., Cherrie,J.B., Mitchell,T.J., Fright,W.R., Ffright,B.C. McCallum,B.C, Evans,T.R.: Reconstruction and representation of 3d objects with radial basis functions. Proc. SIGGRAPH'01, pp. 67–76, 2001.

[4] Duchon,J.: Splines minimizing rotation-invariant semi-norms in Sobolev space. Constructive Theory of Functions of Several Variables, LNCS 571, Springer, 1997.

[5] Hardy,L.R.: Multiquadric equation of topography and other irregular surfaces. Journal of Geophysical Research 76 (8), 1905-1915, 1971.

[6] Fasshauer,G.E.: Meshfree Approximation Methods with MATLAB. World Scientific Publ.2007.

[7] Lazzaro,D., Montefusco,L.B.: Radial Basis functions for multivariate interpolation of large data sets. Journal of Computational and Applied Mathematics, 140, pp. 521-536, 2002.

[8] Macedo,I., Gois,J.P., Velho,L.: Hermite Interpolation of Implicit Surfaces with Radial Basis Functions. Computer Graphics Forum, Vol.30, No.1, pp.27-42, 2011.

[9] Majdisova,Z., Skala,V.: A New Radial Basis Function Approximation with Reproduction, CGVCVIP 2016, pp.215-222, ISBN 978-989-8533-52-4, Portugal, 2016

[10] Majdisova,Z., Skala,V.: A Radial Basis Function Approximation for Large Datasets, SIGRAD 2016, pp.9-14, Sweden, 2016

[11] Nakata,S., Takeda,Y., Fujita,N., Ikuno,S.: Parallel Algorithm for Meshfree Radial Point Interpolation Method on Graphics Hardware. IEEE Trans.on Magnetics, Vol.47, No.5, pp.1206-1209, 2011.

[12] Pan,R., Skala,V.: A two level approach to implicit modeling with compactly supported radial basis functions. Engineering& Comp.,Vol.27.No.3,pp.299-307,ISSN 0177-0667, Springer, 2011.

[13] Pan,R., Skala,V.: Surface Reconstruction with higher-order smoothness. The Visual Computer, Vol.28, No.2, pp.155-162, ISSN 0178-2789, Springer, 2012.

[14] Ohtake,Y., Belyaev,A., Seidel,H.-P.: A multi-scale approach to 3d scattered data interpolation with compactly supported basis functions. In: Proceedings of international conference shape modeling, IEEE Computer Society, Washington, pp 153–161, 2003.

[15] Savchenko,V., Pasco,A., Kunev,O., Kunii,T.L.: Function representation of solids reconstructed from scattered surface points & contours. Computer Graphics Forum, 14(4), pp.181–188, 1995.

[16] Skala,V: Progressive RBF Interpolation, 7th Conference on Computer Graphics, Virtual Reality, Visualisation and Interaction in Africa, Afrigraph 2010, pp.17-20, ACM, 2010

[17] Skala,V., Pan,R.J., Nedved,O.: Simple 3D Surface Reconstruction Using Flatbed Scanner and 3D Print. ACM SIGGRAPH Asia 2013 poster, 2013.

[18] Skala,V.: Projective Geometry, Duality and Precision of Computation in Computer Graphics, Visualization and Games. Tutorial Eurographics 2013, Girona, 2013

[19] Skala,V.: Meshless Interpolations for Computer Graphics, Visualization and Games. Tutorial Eurographics 2015, Zurich, 2015.

[20] Sussmuth,J., Meyer,Q., Greiner,G.: Surface Reconstruction Based on Hierarchical Floating Radial Basis Functions. Computer Graph. Forum, 29(6): 1854-1864, 2010.

[21] Wenland,H.: Scattered Data Approximation. Cambridge University Press, 2010.

[22] Wright,G.B.: Radial Basis Function Interpolation: Numerical and Analytical Developments. University of Colorado, Boulder. PhD Thesis,, 2003

[23] Zapletal,J., Vanecek,P., Skala,V.: RBF-based Image Restoration Utilizing Auxiliary Points. CGI 2009 proceedings, pp.39-44, 2009